\magnification\magstep1
\baselineskip = 18pt
\def\n{\noindent}
 \magnification\magstep1
\baselineskip = 18pt
\def\n{\noindent}
\def\pf{\noindent {\bf Proof.\ \ }}
\overfullrule = 0pt
 \def\qed{{\hfill{\vrule height7pt width7pt
depth0pt}\par\bigskip}}
\def\n{\noindent}
\def \raw{\rightarrow}
\def\cf{{\it cf.\/}\ }
\def\eg{{\it e.g.\/}\ }

\def\cf{{\it cf.\/}\ }

 \overfullrule = 0pt
 \def\qed{{\hfill{\vrule height7pt width7pt
depth0pt}\par\bigskip}}

\def \comp{{\bf C}}

\def\T{{\bf T}}

\overfullrule = 0pt
\def\pf{\medskip\noindent {\bf Proof.}~~}
\def\T{\bf T}
 \def\qed{{\hfill{\vrule height7pt width7pt
depth0pt}\par\bigskip}}

\centerline{\bf The ``maximal" tensor product of operator spaces}
\centerline{by}
\centerline{Timur Oikhberg}

 \centerline{Texas A\&M  University} 

\centerline{College Station, TX 77843, U. S. A.}
\centerline{and}
\centerline{Gilles Pisier\footnote*{Partially supported by the
N.S.F.}}
\centerline{Texas A\&M  University  {and}
 {Universit\'e Paris 6}}

\vskip.5in

\n {\bf Abstract.} In analogy with the  maximal  tensor product of
$C^*$-algebras, we define the ``maximal" tensor product $E_1\otimes_\mu E_2$
of two operator spaces $E_1$ and $E_2$ and we show that it can be
identified completely isometrically with the sum of the two Haagerup
tensor products: \ $E_1\otimes_h E_2 + E_2\otimes_h E_1$. Let $E$ be an
$n$-dimensional operator space. As an application, we show that the
equality $E^* \otimes_\mu E=E^* \otimes_{\rm min} E$ holds isometrically
iff $E = R_n$ or $E=C_n$ (the row or column $n$-dimensional Hilbert
spaces). Moreover, we show that if an operator space $E$ is such that,
for any operator space $F$, we have
$F\otimes_{\min} E=F\otimes_{\mu} E$ isomorphically,
  then $E$ is completely
 isomorphic to either a row or a column  Hilbert  space.

\vfil

\n 1991 AMS Classification Numbers: 47 D 15, 47 D 25, 46 M 05.
\eject

In $C^*$-algebra theory, the minimal and maximal tensor products (denoted
by $A_1 \otimes_{\rm min}~A_2$ and $A_1\otimes_{\rm max} A_2$) of two
$C^*$-algebras $A_1,A_2$, play an important r\^ole, in connection with
``nuclearity'' (a $C^*$-algebra $A_1$ is nuclear if $A_1 \otimes_{\rm min}
A_2 = A_1 \otimes_{\rm max} A_2$ for any $A_2$). See [T] and [Pa1] for more
information and references on this.

In the recently developed theory of operator spaces [ER1-6, Ru1-2, BP1,
B1-3], some specific new versions of the injective and projective
tensor products (going back to Grothendieck for Banach spaces) have been
introduced. 

\n The ``injective'' tensor product of two operator spaces
$E_1,E_2$ coincides with the minimal (or spatial) tensor product and is
denoted by $E_1\otimes_{\rm min} E_2$.

\n Another tensor product of paramount importance for operator spaces is the
Haagerup tensor product,  denoted by
$E_1\otimes_h E_2$ (cf. [CS1-2, PaS, BS]).

\n Assume given two completely isometric embeddings $E_1\subset A_1$,
$E_2\subset A_2$. Then $E_1\otimes_{\rm min}~E_2$ (resp.\ $E_1\otimes_h
E_2$) can be identified with the closure of the algebraic tensor product
$E_1\otimes E_2$ in $A_1\otimes_{\rm min} A_2$ (resp.\ in the ``full''
free product $C^*$-algebra $A_1 * A_2$, see [CES]).
(The ``projective'' case apparently cannot be
described in this fashion and will not be
considered here.)
It is therefore tempting to study the norm induced
on
$E_1\otimes E_2$ by $A_1
\otimes_{\rm max}A_2$. When $A_i = B(H_i)$ $(i=1,2)$ the resulting tensor
product is studied in [JP] and denoted by $E_1\otimes_M E_2$. See also
[Ki] for other tensor products. In the present paper, we follow a
different route:\ we work in the category of (a priori nonself-adjoint)
unital operator algebras, and we use the maximal tensor product in the
latter category (already considered in [PaP]), which extends the
$C^*$-case.

The resulting tensor product, denoted by $E_1\otimes_\mu E_2$ is the
subject of this paper. A brief description of it is as follows:\ we first
introduce the canonical embedding of any operator space $E$ into an
associated ``universal'' unital operator algebra, denoted by $OA(E)$, then
we can define the tensor product $E_1\otimes_\mu E_2$ as the closure of
$E_1\otimes E_2$ in $OA(E_1) \otimes_{\rm max} OA(E_2)$.

Our main result is Theorem 1 which shows that $E_1\otimes_\mu E_2$
coincides with a certain ``symmetrization'' of the Haagerup tensor
product. We apply this (see Corollary~10 and Theorem~16) to find which
spaces $E_1$ have the property that $E_1\otimes_\mu E_2 = E_1 \otimes_{\rm
min} E_2$ for all operator spaces $E_2$.

We refer the reader to the book [Pa1]  for the precise
definitions of all the undefined terminology related to operator spaces
 and complete boundedness, and to [KaR, T] for operator algebras in general. 
 We recall
only that an ``operator space'' is a closed subspace $E\subset B(H)$ of
the $C^*$-algebra of all bounded operators on a Hilbert space $H$. We will
use freely the notion of a completely bounded (in short $c.b.$) map $u\colon
\ E_1\to E_2$ between two operator spaces, as defined e.g.\ in [Pa1]. We
denote by $\|u\|_{cb}$ the corresponding norm and by $cb(E_1,E_2)$ the
Banach space of all c.b.\ maps from $E_1$ to $E_2$.

\n We will denote by $A'$ the commutant of a subset $A\subset B(H)$.

Let $E_1,\ldots, E_n$ be a family of operator spaces. Let
 $\sigma_i\colon \ E_i\to B(H)$
 be complete contractions $(i=1,2,\ldots, n)$. We denote by
$\sigma_1 \cdot\ldots\cdot \sigma_n\colon \ E_1\otimes\cdots \otimes
E_n\to B(H)$ the linear map taking $x_1\otimes\cdots\otimes x_n$ to the
operator $\sigma_1(x_1) \sigma_2(x_2)\ldots \sigma_n(x_n)$.

We define the norm $\|~~\|_\mu$ on $E_1\otimes\cdots\otimes E_n$ as
follows:
$$\|x\|_\mu = \sup\|\sigma_1\cdot\ldots\cdot \sigma_n(x)\|_{B(H)} \leqno 
\forall
x\in E_1\otimes \cdots\otimes E_n$$
where the supremum runs over all possible $H$ and all $n$-tuples
$(\sigma_i)$ of complete contractions as above, with the restriction that
we assume that for any $i\ne j$, the range of $\sigma_i$ commutes with the
range of $\sigma_j$. We will denote by $(E_1\otimes E_2\cdots \otimes
E_n)_\mu$ the completion of $E_1\otimes\cdots\otimes E_n$ for this norm.
In the particular case $n=2$, we denote this simply by $E_1 \otimes_\mu
E_2$, so we have for any $x = \sum x^1_i \otimes x^2_i \in E_1 \otimes E_2$
$$\|x\|_\mu = \sup\{\|\sum \sigma_1 (x^1_i) \sigma_2(x^2_i)\|_{B(H)}\}$$
where the supremum runs over all possible pairs $(\sigma_1,\sigma_2)$ of
complete contractions (into some common $B(H)$) with commuting ranges,
i.e.\ such that $\sigma_1(x_1) \sigma_2(x_2) = \sigma_2(x_2)
\sigma_1(x_1)$ for all $x_1\in E_1$, $x_2\in E_2$.

The space $(E_1\otimes \cdots \otimes E_n)_\mu$ can obviously be equipped
with an operator space structure associated to the embedding
$$J\colon \ (E_1\otimes \cdots \otimes E_n)_\mu \to \bigoplus_\sigma
B(H_\sigma)\subset B(\bigoplus_\sigma H_\sigma)$$
where the direct sum runs over all $n$-tuples $\sigma = (\sigma_1,\ldots,
\sigma_n)$ with $\sigma_i\colon \ E_i\to B(H_\sigma)$ such that
$\|\sigma_i\|_{cb}\le 1$ and the $\sigma_i$'s have commuting ranges.

We note in passing that if we define 
$\hat \sigma_i\colon \ E_i\to \bigoplus_\sigma
B(H_\sigma)\subset B(\bigoplus_\sigma H_\sigma)$
by
$\hat \sigma_i(x)=  \oplus_\sigma \sigma_i(x)$, then we have
$J(x)= \hat \sigma_1 ... \hat \sigma_n$ and the maps $\hat \sigma_i$ have commuting ranges.

Thus we can now unambiguously refer to $(E_1\otimes \cdots \otimes
E_n)_\mu$, and in particular to $E_1\otimes_\mu E_2$ as operator spaces.

We will give more background on operator spaces and c.b.\ maps below.
For the moment, we merely define a ``complete metric surjection'':\ by
this we mean a surjective mapping $Q\colon \ E_1\to E_2$ between
two operator spaces, which induces a complete isometry from $E_1/\ker(Q)$
onto $E_2$.

To state our main result, we also need the notion of $\ell_1$-direct sum of
two operator spaces $E_1,E_2$:\ this is an operator space denoted by $E_1
\oplus _1 E_2$. The norm on the latter space is as in the usual
$\ell_1$-direct sum, i.e.\ we have
$$\|(e_1,e_2)\|_{E_1\oplus_1E_2} = \|e_1\|_{E_1} + \|e_2\|_{E_2}$$
but the operator space structure is such that for any pair $u_1\colon \
E_1\to B(H)$, $u_2:\ E_2\to B(H)$ of complete contractions, the mapping
$(e_1,e_2) \to u_1(e_1) + u_2(e_2)$ is a complete contraction from $E_1
\oplus_1 E_2$ to $B(H)$.

The simplest way to realize this operator space $E_1\oplus_1E_2$ as a
subspace of $B({\cal H})$ for some ${\cal H}$ is to consider the
collection $I$ of all pairs $p=(u_1,u_2)$ as above with $H=H_p$ (say) and
to define the embedding
$$J\colon \ E_1\oplus_1 E_2 \to B\left( \bigoplus_{p\in I} H_p\right)$$
defined by $J(e_1,e_2) = \bigoplus\limits_{(u_1,u_2)\in I} [u_1(e_1) +
u_2(e_2)]$. Then, we may as well define the operator space structure of
$E_1 \oplus_1 E_2$ as the one induced by the isometric embedding $J$. In
other words, $E_1\oplus_1 E_2$ can be viewed as the ``maximal'' direct sum
for operator spaces, in accordance with the general theme of this paper.

We will denote by $E_1\otimes E_2$ the linear tensor
product of two vector spaces and by $v \to {}^tv$ the transposition map,
i.e.\ for any $v = \sum x_i\otimes y_i$
in $E_1\otimes E_2$, we set ${}^t v =
\sum y_i\otimes x_i$. The identity map on a space $E$ will be
denoted by $Id_E$.

\n We will denote by $E_1\otimes_h E_2$ the Haagerup tensor product of
two operator spaces for which we refer to [CS1-2, PaS].  

\n {\bf Convention:}\ We reserve the term ``{\it morphism\/}'' for a
unital completely contractive homomorphism $u\colon \ A\to B$ between two
unital operator algebras.

Our main result is the following one.

\proclaim Theorem 1. Let $E_1$, $E_2$ be two operator spaces.
Consider the mapping
$$Q\colon \ (E_1\otimes_h E_2) \oplus_1 (E_2 \otimes_h E_1) \to
E_1\otimes_\mu E_2$$
defined on the direct sum of the linear tensor products
 by $Q(u\oplus v) = u+{}^tv$. Then $Q$ extends to a complete metric
surjection from $(E_1\otimes _h E_2) \oplus_1 (E_2\otimes_hE_1)$ onto
$E_1 \otimes_\mu E_2$. In particular, for any $u$ in $E_1\otimes E_2$, we
have:\ $\|u\|_\mu<1$ iff there are $v,w$ in $E_1\otimes E_2$ such that
$u = v+w$
and
$\|v\|_{E_1\otimes_h E_2} + \|{}^tw\|_{E_2\otimes_hE_1} < 1.$

\n {\bf Remark.} In the terminology of [P2], the preceding statement means
that $E_1\otimes_\mu E_2$ is completely isometric to the ``sum'' (in the
style of interpolation theory, see [P2]) $E_1\otimes_h E_2 + E_2\otimes_h
E_1$ (in analogy with $R+C$).

\n {\bf Remark.} We first recall a simple consequence of the Cauchy-Schwarz inequality:\ for any
$a_1,\ldots, a_n$, $b_1,\ldots, b_n$ in a $C^*$-algebra $A$, we have
$$\big\|\sum a_ib_i\big\| \le \big\|\sum a_ia^*_i\big\|^{1/2} \big\|\sum
b^*_ib_i\big\|^{1/2}.$$
Hence if $a_ib_i = b_ia_i$ we also have
$$\big\|\sum a_ib_i\big\| \le \big\|\sum a^*_ia_i\big\|^{1/2} \big\|\sum
b_ib^*_i\big\|^{1/2}.$$
From these it is easy to deduce that the above map $Q$ is  completely contractive.

Using [CS1-2, PaS], one can see that the following  statement is a dual reformulation:

\proclaim Theorem 2. Let $E_1,E_2$ be two operator spaces,
 and let
$\varphi\colon \ E_1\otimes E_2 \to B( {\cal H})$ be a linear mapping. The 
following
are equivalent:
\item{\rm (i)} $\|\varphi\|_{cb(E_1\otimes_\mu E_2, B({\cal H})) } \le 1$.
\item{\rm (ii)} For some Hilbert space ${ H}$, there are complete contractions
$\alpha_1\colon \ E_1\to B(H,{\cal H})$,
$\alpha_2\colon \ E_2\to B({\cal H},H)$,
 and $\beta_1\colon \ E_1 \to B({\cal H},H)$,
$\beta_2\colon \ E_2 \to B(H,{\cal H})$,
such that
$$\varphi(x_1\otimes x_2) = \alpha_1(x_1)
\alpha_2(x_2) = \beta_2 (x_2)
\beta_1(x_1).\leqno \forall (x_1,x_2)\in E_1\times E_2$$
\item{\rm (iii)} For some Hilbert space ${ H}$, there are complete contractions
$\sigma_i\colon
\ E_i\to B(H)$   $(i=1,2)$, with commuting ranges, and contractions
$V\colon\ {\cal H}\to H$ and $W\colon\ {  H}\to {\cal H}$, such that
$$\varphi(x_1\otimes x_2) = W\sigma_1(x_1)
\sigma_2(x_2)V.\leqno \forall (x_1,x_2)\in E_1\times E_2$$

\pf Assume (i). By the preceding remark, $\varphi$ defines a complete contraction into $B({\cal
H})$ both from
$E_1\otimes_h E_2$ 
and from $E_2\otimes_h E_1$. Then (ii) follows from the Christensen-Sinclair factorization
theorem for bilinear maps,
  extended to general operator spaces by Paulsen and Smith in [PaS].

\n Now assume (ii).  Let $H_1 = {\cal H}$, $H_2=H$ and $H_3={\cal H}$. We define
maps $\sigma_1\colon \ E_1 \longmapsto B(H_1\oplus H_2 \oplus H_3)$ and
$\sigma_2\colon \ E_2 \longrightarrow B(H_1\oplus H_2 \oplus H_3)$ using
matrix notation, as follows
$$\eqalign{\sigma_1(x_1) &= \left(\matrix{0&\alpha_1(x_1)&0\cr
0&0&\beta_1(x_1)\cr 0&0&0\cr}\right)\cr
\sigma_2(x_2) &= \left(\matrix{0&\beta_2(x_2)&0\cr
0&0&\alpha_2(x_2)\cr 0&0&0\cr}\right).}$$
Then, by (ii) we have
$$\sigma_1(x_1) \sigma_2(x_2) = \sigma_2(x_2) \sigma_1(x_1) =
\left(\matrix{0&0&\varphi(x_1\otimes x_2)\cr
0&0&0\cr 0&0&0\cr}\right)$$
hence $\sigma_1,\sigma_2$ have commuting ranges, and are complete
contractions.

\n Therefore if we let $W\colon \ H_1\oplus H_2\oplus H_3\to {\cal H}$ be the
projection onto the first coordinate and $V\colon \ {\cal H} \to H_1
\oplus H_2\oplus H_3$ be the isometric inclusion into the third
coordinate, then we obtain (iii).

\n Finally, the implication (iii) implies (i) is obvious by the very definition of $E_1\otimes_\mu
E_2$. \qed

 \n{\bf First proof of Theorem 1.} By duality, it clearly suffices to show that 
for any linear map $\varphi\colon \ E_1\otimes E_2 \to B( {\cal H})$
the norms   $\|\varphi\|_{cb(E_1\otimes_\mu E_2 \to B( {\cal H})) }$ and
$\|\varphi Q\|_{cb}$
are equal. But this is precisely the meaning of the equivalence
between (i) and (ii) in Theorem 2. Thus we conclude that Theorem 2 implies Theorem 1. \qed   

The main idea of the second proof of Theorem 1 is to use the universal unital
operator algebras of operator spaces as initiated in  [P4], to
relate their free product with their ``maximal'' tensor product, and to
use the appearance of the Haagerup tensor product inside the free product. While that proof is
longer,  the steps are of independent interest and somehow the principle behind it, which 
  can be applied in   other instances,   is perhaps easier to generalize than
the sort of ``trick" used in the  proof of Theorem 2. For that reason, we find it worthwhile to
include it.

Let $E$ be an operator space. Let $T(E) = \comp \oplus E\oplus (E\otimes
E) \oplus\cdots$ be its tensor algebra, so that any $x$
 in $T(E)$ is a sum $x = \sum x_n$ with $x_n \in
E\otimes\cdots \otimes E$ ($n$ times) with $x_n=0$ for all $n$
sufficiently large.
For each linear $\sigma\colon\ E\to B(H)$ we denote by $T(\sigma)\colon \
T(E)\to B(H)$ the unique unital homomorphism extending $\sigma$.

Let $C$ be the collection of all $\sigma\colon \ E\to B(H_\sigma)$ with
$\|\sigma\|_{cb}\le 1$. We define an embedding
$$J\colon \ T(E) \to B\left(\bigoplus_{\sigma\in C} H_\sigma\right)$$
by setting
$$J(x) = \bigoplus_{\sigma\in C} T(\sigma)(x).$$
Then $J$ is a unital homomorphism. We denote by $OA(E)$ the unital
operator algebra obtained by completing $J(T(E))$. We will always view
$T(E)$ as a subset of $OA(E)$, so we identify $x$ and $J(x)$ when $x\in
T(E)$. Observe that the natural inclusion
$$E\to OA(E)$$
is obviously a complete isometry. More generally, the natural inclusion of
$E\otimes\cdots\otimes E$ ($n$ times) into $OA(E)$ defines a completely
isometric embedding of $E \otimes_h\cdots \otimes_h E$ into $OA(E)$. (This
follows from a trick due to Varopoulos, and used by Blecher in [B1], see
[P4] for details).

The algebra $OA(E)$ is characterized by the following universal
property:\ for any $\sigma\colon \ E\to B(H)$ with $\|\sigma\|_{cb}\le 1$,
there is a unique morphism $\hat\sigma\colon \ OA(E) \to B(H)$ extending
$\sigma$ (here we view $E$ as embedded into $OA(E)$ in the natural way). See [Pes] for the
self-adjoint analogue.
 
We now turn to the free product in the category of unital operator
algebras. Let $A_1,A_2$ be two such algebras and let ${\cal F}$ be their
algebraic free product as unital algebras (i.e.\ we identify the units and
``amalgamate over $\comp$''). For any pair $u = (u_1,u_2)$ of morphisms,
as follows $u_i\colon \ A_i\to B(H_u)$ $(i=1,2)$, we denote by
$u_1*u_2\colon \ {\cal F}\to B(H_u)$ the unital homomorphism extending
$u_1,u_2$ to the free product.

\n Then we consider the embedding
$$J\colon \ {\cal F} \to \bigoplus _u B(H_u)$$
defined by $J(x) = \bigoplus\limits_u [u_1*u_2(x)]$, for all $x$ in ${\cal
F}$. Note that $J$ is a unital homomorphism.

We define the free product $A_1*A_2$ (in the category of unital operator algebras) as the
closure of $J({\cal F})$. Actually, we will identify ${\cal F}$ with
$J({\cal F})$ and consider that $A_1*A_2$ is the completion of ${\cal F}$
relative to the norm induced by $J$. Moreover, we will consider $A_1*A_2$
as a unital operator algebra equipped with the operator space structure induced by $J$.

It is easy to see that $A_1*A_2$ is characterized by the (universal)
property that for any pair of morphisms
$u_i\colon \ A_i\to B(H)$ $(i=1,2)$ there is a unique morphism from
$A_1*A_2$ to $B(H)$ which extends both $u_1$ and $u_2$.

We now turn to the maximal tensor product in the category of unital
operator algebras. This is defined in [PaP], so we only briefly recall the
definition:\ Let $A_1,A_2$ be two unital operator algebras. For any pair
$\pi = (\pi_1,\pi_2)$ of morphisms $\pi_i\colon \ A_i \to B(H_\pi)$ with
commuting ranges, we denote by $\pi_1\cdot\pi_2$ the morphism from
$A_1\otimes A_2$ to $B(H_\pi)$ which takes $a_1\otimes a_2$ to
$\pi_1(a_1)\pi_2(a_2)$. Then we consider the embedding
$$J\colon \ A_1\otimes A_2 \to B\left(\bigoplus_\pi H_\pi\right)$$
defined by $J(x) = \bigoplus\limits_\pi \pi_1\cdot\pi_2(x)$.  We define
$$\|x\|_{\rm max} = \sup_\pi \|\pi_1\cdot \pi_2(x)\|$$
and we denote by $A_1\otimes_{\rm max}A_2$ the completion of $A_1\otimes
A_2$ for this norm. We will consider $A_1\otimes_{\rm max} A_2$ as a
unital operator algebra, using the isometric embedding $J$ just defined.

We will use several elementary facts which essentially all follow from the
universal properties of the objects we have introduced.

\proclaim Lemma 3. $OA(E_1) * OA(E_2) \simeq OA(E_1\oplus_1 E_2)$
completely isometrically.

\n {\bf Remark.} The ``functor'' $E\to OA(E)$ is both injective and
projective:\ indeed, if $E_2~\subset~E_1$ is a closed subspace then the
associated maps $OA(E_2) \to OA(E_1)$ and $OA(E_1) \to\break OA(E_1/E_2)$ are
respectively a complete isometry and a complete metric surjection.
\medskip

\proclaim Lemma 4. Let $A_1,A_2$ be two unital operator algebras. Then the
natural morphism $Q\colon \ A_1*A_2\to A_1 \otimes_{\rm max} A_2$ is a
complete metric surjection. More precisely, the restriction of $Q$ to the
algebraic free product ${\cal F}$ defines a complete isometry between
${\cal F}/\ker(Q_{|{\cal F}})$ and $A_1\otimes A_2 \subset A_1
\otimes_{\rm max} A_2$.

\pf  Let $j_i\colon \ A_i\to {\cal F}$ $(i=1,2)$
denote the canonical inclusion maps. By the universal property of the free
product, there is a unique morphism $Q$ extending the inclusion mappings
$$a_1\to a_1\otimes 1\quad \hbox{and}\quad a_2\to 1\otimes a_2.$$
Now consider $Q_{|{\cal F}}\colon \ {\cal F}\to A_1\otimes_{\rm max} A_2$
and the associated injective morphism 
$$\widehat Q\colon\ {\cal F}/\hbox{ker}(Q_{|{\cal F}}) \to
A_1\otimes_{\rm
max} A_2.$$ 
Let $q\colon \ {\cal F}\to {\cal F}/\hbox{ker}(Q_{|{\cal F}})$
be the canonical surjection so that $Q_{|{\cal F}} = \widehat Qq$.
Clearly, the range of $\widehat Q$ coincides with $A_1\otimes A_2$.
Let $\sigma\colon \ A_1\otimes A_2 \to {\cal F}/\hbox{ker}(Q_{|{\cal F}})$
be the restriction of $\widehat Q^{-1}$ to $A_1\otimes A_2$. Then $\sigma$
is a unital homomorphism into an operator algebra (by [BRS]) such that
$a_1\to \sigma(a_1\otimes 1)$ and $a_2\to \sigma(1\otimes a_2)$ are
morphisms on $A_1$ and $A_2$ respectively. (Indeed, $a_1\to
\sigma(a_1\otimes 1)$ coincides with the composition
$A_1 {\buildrel j_1\over \longrightarrow} {\cal
F} {\buildrel q\over
\longrightarrow} {\cal F}/\ker(Q_{|{\cal F}}),$
and similarly for $A_2$.)
Therefore, by definition of $A_1\otimes_{\rm max}
A_2$, $\sigma$ must extend to a morphism from
$A_1\otimes_{\rm max} A_2$ to ${\cal
F}/\ker(Q_{|{\cal F}})$. In particular, this
implies that $\widehat Q$ is a complete isometry
from ${\cal F}/\ker(Q_{|{\cal F}})$ onto
$A_1\otimes A_2$ viewed as an operator subspace
of $A_1\otimes_{\rm max} A_2$. A fortiori, $Q$ is
a complete metric surjection onto
$A_1\otimes_{\rm max} A_2$.\qed\medskip

\proclaim Lemma 5. The natural inclusion of $E_1\otimes_\mu E_2$ into
$OA(E_1) \otimes_{\rm max} OA(E_2)$ is a completely isometric embedding.

The next lemma (already used in [B1]) is elementary.

\proclaim Lemma 6. Consider an element
$$x = x_0 + x_1 +\cdots+x_n+\cdots$$
in $T(E)$, with $x_n \in E\otimes\cdots \otimes E$ ($n$ times). Then the
mapping $x\to x_n$ defines a completely contractive projection on $OA(E)$.

\pf Let $m$ denote the normalized Haar measure on the unidimensional torus
$\T$. For $z$  in $\T$, let $x(z) = \sum\limits_{n\ge 0} z^nx_n$. By
definition of $OA(E)$, we clearly have
$$\|x(z)\| = \|x\|,$$
hence
$$\|x_n\| = \left\|\int z^{-n} x(z) m(dz)\right\|\le \|x\|.$$
This shows that $x\to x_n$ is a contractive linear projection. The
argument for complete contractivity is analogous and left to the
reader.\qed

The next result, which plays an important r\^ole in the sequel, might be
of
independent interest.

\proclaim Lemma 7. Let $E_1,E_2$ be two operator spaces. Let $X = (E_1
\oplus_1 E_2) \otimes_h (E_1\oplus_1 E_2)$. With the obvious
identifications, we may view $E_1\otimes E_2 + E_2\otimes E_1$ as a linear
subspace of $X$. Let $S$ be its closure in $X$. Then we have
$$S \simeq (E_1 \otimes_h E_2) \oplus_1 (E_2\otimes_h E_1)$$
completely isometrically.

\pf Obviously we have completely contractive natural inclusions $E_1
\otimes_h E_2\to X$ and $E_2\otimes_h E_1\to X$, whence a natural inclusion
$(E_1\otimes _h E_2) \oplus_1(E_2\otimes_h E_1) \to X$. To show that this
is completely isometric it clearly suffices to show that $S$ has the
``universal'' property characteristic of the $\oplus_1$-direct sum.
Equivalently, it suffices to show that every completely contractive
mapping $\sigma\colon \ (E_1\otimes_h E_2) \oplus_1 (E_2 \otimes_h E_1)
\to B(H)$
defines a completely contractive mapping from $S$ to $B(H)$ (then we may
apply this when $\sigma$ is a completely isometric embedding). So let
$\sigma$ be such a map. Clearly, we can assume that $\sigma(x\otimes y) =
u(x) + v(y)$ with $u\colon \ E_1\otimes_h E_2\to B(H)$, $v\colon \
E_2\otimes_h E_1\to B(H)$ such that $\|u\|_{cb} \le 1$, $\|v\|_{cb}\le 1$.
By the factorization of $c.b.$-bilinear maps ([CS, PaS]) we can further write
$u(x_1\otimes x_2) = u_1(x_1) u_2(x_2)$ and $v(y_2\otimes y_1) = v_2(y_2)
v_1(y_1)$ where $u_i \colon\ E_i\to B(H)$ and $v_i\colon \ E_i\to B(H)$ are
all completely contractive.
Let us then define $\alpha\colon \ E_1\oplus_1E_2\to B(H)$ and
$\beta\colon \ E_1\oplus_1 E_2\to B(H)$ by\break $\alpha(x_1\oplus x_2)
=u_1(x_1) + v_2(x_2)$ and $\beta(x_1\oplus x_2)= v_1(x_1) + u_2(x_2)$.  
By
definition of $E_1\oplus_1E_2$, these maps are still complete
contractions. Moreover, we have for any $z$ in $S$, say $z = x+y$ with
$x\in E_1\otimes E_2$, $y\in E_2\otimes E_1$
$$\alpha\cdot\beta(z) = u(x) + v(y) = \sigma(z).$$
hence we conclude that $\sigma$ admits an extension $\tilde\sigma$ (namely
$\tilde\sigma = \alpha\cdot\beta$) defined on the whole of $X$ with
$$\|\tilde\sigma\|_{cb(X,B(H))} \le \|\alpha\|_{cb} \|\beta\|_{cb} \le
1,$$
a fortiori $\|\sigma\|_{cb(S,B(H))} \le 1$.\qed

\proclaim Lemma 8. Let $A_1,\ldots, A_n$ be unital operator algebras. Then
the natural mapping from $A_1\otimes A_2\cdots \otimes A_n$ to
$A_1 * A_2 *\cdots * A_n$ defines a completely isometric embedding of
 $A_1\otimes_h
A_2 \cdots  \otimes_h A_n$ into $A_1 * A_2 *\cdots * A_n$.

\pf In essence, this is proved in [CES], but only for the non-unital free
product. The unital case is done in detail in [P3] so we skip it.\qed

\n {\bf Second proof of Theorem 1.}  Consider $u$ in $E_1 \otimes  E_2$ with
$\|u\|_\mu<1$. Let ${\cal F}$ be, as before, the algebraic free product of
 $OA(E_1)$ and $OA(E_2)$. By Lemma~5, we have $\|u\|_{OA(E_1)\otimes_{\rm max}
OA(E_2)}<1$, hence by Lemma~4, there is an element $\hat u$ in ${\cal F}$
 with $\|\hat u\| <1$ such that $Q(\hat u) = u$. By Lemma~3 we may
write as well
$$\|\hat u\| _{OA(E_1\oplus_1 E_2)} <1.$$
Let us write $\hat u = \hat u_0 + \hat u_1 + \hat u_2 +\cdots$ where $\hat
u_d \in (E_1\oplus_1 E_2) \otimes_h \cdots \otimes_h (E_1 \oplus_1 E_2)$
($d$ times). By Lemma~6 we have $\|\hat u_d\|<1$.

\n Let $z_1,z_2$ be complex numbers with $|z_i|\le 1$. There is a unique
morphism $\pi_{z_i}\colon\ OA(E_i) \to OA(E_i)$ extending $z_i Id_{E_i}$.
We will use the morphisms
$$\pi_{z_1} \otimes \pi_{z_2}\quad \hbox{acting on}\quad OA(E_1)\otimes_
{\rm max}
OA(E_2)$$
and
$$\pi_{z_1} * \pi_{z_2} \quad \hbox{acting on}\quad OA(E_1) * OA(E_2).$$
Note that we trivially have the following relation:
$$\pi_{z_1} \otimes \pi_{z_2} = Q\circ [\pi_{z_1}* \pi_{z_2}].$$
It follows that
$$z_1z_2u = Q[\pi_{z_1}* \pi_{z_2}(\hat u)].$$
Hence identifying the coefficient of $z_1z_2$ on the right hand side we
obtain
$$u = Q[\tilde u]$$
where $\tilde u$ is in the subspace
$$S = {\rm span}[E_1\otimes E_2 + E_2\otimes E_1] \subset
(E_1\oplus_1 E_2) \otimes_h (E_1\oplus_1E_2)$$
considered in Lemma~7, and where we view $(E_1\oplus_1 E_2) \otimes_h
(E_1\oplus_1E_2)$ as the subspace of $OA(E_1\oplus_1 E_2)$ formed of all terms
of degree $2$, according to Lemma 8.

\n Hence by Lemma~7, we conclude that $\tilde u$ can be written as $v+w$ with
$v\in {E_1\otimes E_2}$ and $w\in {E_2\otimes E_1}$ such that  $$\|v\|_{E_1\otimes_h E_2} +
\|w\|_{E_2\otimes_h E_1} < 1.$$ This shows that $\forall u\in E_1\otimes  E_2$ with $\|u\|_\mu <1$
there are $v, w$ as above with $u=v+{}^t w$.

\n Thus the natural mapping is a metric surjection from $E_1\otimes_h E_2
\oplus_1 E_2 \otimes_h E_1$ onto $E_1\otimes_\mu E_2$. To show that this is
a complete surjection, one simply repeats the argument with
$M_n(E_1\otimes_\mu E_2)$ instead of $E_1\otimes_\mu E_2$. We leave the
easy details to the reader.\qed

\n {\bf Remark.} Following [Di], we say that a collection of  Banach  
algebras which is stable by arbitrary direct sums, subalgebras and quotients is a variety. Let
${\cal V}$ be a variety formed
of unital operator algebras. Of course, we are interested in their operator space  (and not only
their Banach space) structure and we use unital completely contractive homomorphisms as morphisms.
One of the advantages of the second proof over the first one is that its ``categorical principle"
can be easily adapted to compute the analogue of the
$\mu$-tensor product obtained when one restricts all maps to take their values into an algebra
belonging to some fixed given variety ${\cal V}$.

\n {\bf Remark.} Using the same techniques, one can prove the {\it isomorphic}
version of Theorem 1 for more than two spaces: namely, that
$(E_1 \otimes E_2 \otimes \ldots \otimes E_n)_\mu$ is
$f(n)$-completely isomorphic to $ \sum_{\sigma \, {\rm permutation}}
E_{\sigma(1)} \otimes_h \ldots \otimes_h E_{\sigma(n)}$. We do not believe
that the {\it isometric} analog of Theorem~1 holds true for $n>2$. Nevertheless,
generalizing Lemma~7, we can prove that, if
$$
X = (E_{11} \oplus_1 \ldots \oplus_1 E_{1n}) \otimes_h \ldots \otimes_h
(E_{K1} \oplus_1 \ldots \oplus_1 E_{kn}) ;
$$
$$
Y = E_{11} \otimes_h \ldots \otimes_h E_{k1} \oplus_1  \ldots \oplus_1
E_{1n} \otimes_h \ldots \otimes_h E_{kn} ;
$$
$$
S = E_{11} \otimes_h \ldots \otimes_h E_{k1} \oplus \ldots \oplus
E_{1n} \otimes_h \ldots \otimes_h E_{kn} \hookrightarrow X ,
$$
then $Y \simeq S$ completely isometrically.
However, to emulate the second  proof of Theorem~1,
we would need a stronger version of Lemma~7, as follows:\hfill \break
if
$
X = (E_1 \oplus_1 \ldots \oplus_1 E_n) \otimes_h \ldots \otimes_h
(E_1 \oplus_1 \ldots \oplus_1 E_n) \, \, \, ( n \, {\rm times} )$ ;
$Y = \oplus_1 E_{\sigma(1)} \otimes_h \ldots \otimes_h E_{\sigma(n)}$
(the sum being taken over all permutations $\sigma$),
 and \hfill \break $
S = \sum_{\sigma \, {\rm permutation}}
E_{\sigma(1)} \otimes \ldots \otimes E_{\sigma(n)}
\hookrightarrow X ,$
then
$$ S \simeq Y \quad \hbox{(completely isometrically).}
\leqno(1)$$

\n However, this is not always true.
Indeed, consider $n=2k$, $N=k!$, $E_1 = \ldots = E_n = \comp$; let
$e_1, \ldots, e_n$ be  unit vectors in $E_1, \ldots, E_n$. If $\tau$ and
$\pi$ are permutations of $[1, \ldots, k]$ and $[k+1, \ldots, n]$,
respectively, let $x_\tau = e_{\tau(1)} \otimes \ldots \otimes e_{\tau(k)}$,
$y_\pi = e_{\pi(k+1)} \otimes \ldots \otimes e_{\pi(n)}$. In this case,
$X = \ell_1^n \otimes_h \ldots \otimes_h \ell_1^n$, $Y = \ell_1^{n!}$,
and $Z={\rm span}[x_\tau \otimes y_\pi \, | \, \tau, \pi {\rm \, permutations}]
\hookrightarrow X$ is a subspace of $S$.
To contradict (1), we will show that the natural
identity map $id : Z \raw \ell_1^{N^2}$ is not a complete isometry.
Indeed, if it were, ${\rm span}[x_\tau]$ and ${\rm span}[y_\pi]$ would
be completely isometric to $\ell_1^N$; hence,
$|| id : \ell_1^N \otimes_h \ell_1^N \raw \ell_1^{N^2} || \leq 1$, or, by
duality,
$$
|| id : \ell_\infty^{N^2} \raw \ell_\infty^N \otimes_h \ell_\infty^N || =1.
\leqno{(2)}$$
Note that $\ell_\infty^{N^2}$ can be canonically identified with
$B(\ell_1^N, \ell_\infty^N)$ and, by Theorem~3.1 of [B1],
$\ell_\infty^N \otimes_h \ell_\infty^N =
\Gamma_2 (\ell_1^N, \ell_\infty^N)$. Therefore, (2) implies that
$B(\ell_1^N, \ell_\infty^N) = \Gamma_2 (\ell_1^N, \ell_\infty^N)$, which
is not true (\cf\eg [P1, p. 48]). Hence, (1) is false, and
this `multispace' version of Lemma~7 does not hold.

\n {\bf Remark.}
If $X,Y$ are Banach spaces, and if $v\in Y\otimes X$, let us
denote by
$\gamma_2(v)$ the norm of factorization through Hilbert space of the linear map $\tilde v\colon \
Y^*\to X$ associated to $v$. This is  a classical notion in Banach space theory (\cf\eg [P1, p.
21]). Note that Theorem 1 obviously implies that for any $v$ in  $E_1 \otimes E_2$ ($E_1, E_2$
being arbitrary operator spaces), we have $$\gamma_2(v)\le \|v\|_\mu.\leqno(3)$$

\n {\bf Remark.}
Note that
$$
(E_1 \otimes_\mu E_2) \otimes_\mu E_3 =
E_1 \otimes_h E_2 \otimes_h E_3 +
E_3 \otimes_h E_1 \otimes_h E_2 +
E_2 \otimes_h E_1 \otimes_h E_3 +
E_3 \otimes_h E_2 \otimes_h E_1 ,
$$
but the above expression does not necessarily coincide with
$(E_1 \otimes E_2 \otimes E_3)_\mu$, and moreover the $\mu$-tensor product is not associative,
in sharp contrast with the Haagerup one (or with the maximal tensor product for unital operator
algebras). In particular, in general the natural mapping from $E_1 \otimes_\mu (E_2 \otimes_\mu
E_3)$ into
$(E_1 \otimes_\mu E_2) \otimes_\mu E_3$  is unbounded  (and actually only makes sense   on the
linear tensor products). All this follows from the   counterexample below, kindly communicated to
us by C. Le Merdy. Let
$X$ be a Banach space and let $K$ be the algebra of all compact operators on $\ell_2$.  Take
$E_1=C$,
$E_2=R$, and
$E_3=\min(X)$ in the sense of [BP1]. Assume that  we have a bounded map
$$(E_1 \otimes E_2 \otimes E_3)_\mu \to (E_1 \otimes_\mu E_2) \otimes_\mu E_3 .$$
Then a fortori we have a bounded map
  $E_1 \otimes_\mu (E_2 \otimes_\mu E_3)\to (E_1 \otimes_\mu
E_2) \otimes_\mu E_3$, and consequently a bounded map
 $E_1 \otimes_h E_3 \otimes_h E_2 \to (E_1 \otimes_\mu E_2)
\otimes_\mu E_3 .$ But then $C \otimes_h E_3 \otimes_h R=K\otimes_{\min} E_3$ 
completely isometrically (see [BP1, ER4]), hence it is isometric to the (Banach space theoretic)
injective tensor product $K\check{\otimes} X$. Moreover, since $R\otimes_h C$ is
isometric to
$K^*$, by Theorem 1 we have
$C\otimes_\mu R=K$ isometrically.

\n  Thus,   we would have a bounded map
from  $K\check{\otimes} X$ to $(C\otimes_\mu R)\otimes_\mu E_3$, and this  would imply by (3),
that for some constant $C$,  for all $v$ in $K\check{\otimes} X$,
 we would have $\gamma_2(v)\le C\|v\|_{\vee}$. However, it is well known that
this fails  at least for some Banach space $X$ (take   for example $X=\ell_1$ and $v=\sum_1^n
e_{ii}\otimes e_i$, so that $v$ represents an isometric embedding of
$\ell_n^\infty$ into $K$, then $\|v\|_{\vee}=1$ and $\gamma_2(v)=\sqrt{n}$, \cf [P1, p. 48] for
more on this question).

\vskip .1 in

We now give several consequences and reinterpret Theorem~1, in terms of
factorization.

The following notation will be convenient. Let $X$ be an operator space.
We will say that a linear map $u\colon \ E_1\to E_2$ between operator
spaces factors through $X$ if there are maps $w\colon \ E_1\to X$ and
$v\colon \ X\to E_2$ such that $u=vw$. We will denote by $\Gamma_X(E_1,E_2)$
the class of all such mappings and moreover we let
$$\gamma_X (u) = \inf\{\|v\|_{cb} \|w\|_{cb}\}$$
where the infimum runs over all possible such factorizations. Let us
denote by ${\cal K}$ the $C^*$-algebra of all compact operators on
$\ell_2$, with its natural ``basis'' $(e_{ij})$.

The preceding notation applies in particular when $X = {\cal K}$ and gives
us the space $\Gamma_{\cal K}(E_1,E_2)$. In the case $X={\cal K}$, it is
easy to check that $\gamma_{\cal K}$ is a norm with which $\Gamma_{\cal
K}(E_1,E_2)$ becomes a Banach space.

We wish to relate the possible factorizations of a map through ${\cal K}$
with its possible factorizations through two specific subspaces of ${\cal
K}$, namely the row and column Hilbert spaces defined by
$$\eqalign{R &= \overline{\rm span}(e_{1j}\mid j=1,2,\ldots)\cr
C &= \overline{\rm span}(e_{i1}\mid i=1,2,\ldots).}$$
Clearly these subspaces of ${\cal K}$ admit a natural completely
contractive projection onto them (namely $x\to e_{11}x$ is a projection
onto $R$, and $x\to xe_{11}$ one onto $C$). Therefore we have
$\gamma_{\cal K}(Id_R) = 1$ and
$\gamma_{\cal K}(Id_C) = 1$. A fortiori any linear map $u\colon \ E_1\to
E_2$ which factors either through $R$ or through $C$ factors through
${\cal K}$ and we have
$$\gamma_{\cal K}(u) \le \gamma_R(u)\quad \hbox{and}\quad \gamma_{\cal
K}(u) \le \gamma_C(u).$$
A fortiori, if $u=v+w$ for some $v\colon\ E_1\to E_2$ and $w\colon \
E_1\to E_2$, we have
$$\gamma_{\cal K}(u) \le \gamma_R(v) + \gamma_C(w).$$
Note that if $v$ and $w$ are of finite rank, then with the obvious
identifications, we have
$$\gamma_R(v)=\|v\|_{E_1^*\otimes_h E_2}\quad \hbox{and} \quad
\gamma_C(w)=\|{}^t w\|_{E_2\otimes_h E_1^*} .$$ Thus, from Theorem 1 we deduce:

\proclaim Corollary 9. Let $E_1,E_2$ be operator spaces. Consider $u$ in
$E^*_1\otimes E_2$ and let $\tilde u\colon \ E_1\to E_2$ be the associated
finite rank operator. Then we have
$$\gamma_{\cal K}(\tilde u) \le \|u\|_\mu \leqno(4)$$

\proclaim Corollary 10. Let $E$ be an $n$-dimensional operator space. Let
$i_E \in E^*\otimes E$ be associated to the identity of $E$
and let
$$\mu(E) = \|i_E\|_\mu.$$
Then
$$\max\{\gamma_{\cal K}(Id_E), \gamma_{\cal K}(Id_{E^*})\} \le \mu(E).
\leqno  {  (5)}$$
Moreover $\mu(E)=1$ iff either $E = R_n$ or $E=C_n$ (completely
isometrically).

\pf Note that (5) clearly follows from (4). Assume that $\mu(E)=1$.
Then by Theorem~1
(and an obvious compactness argument) we have a decomposition $Id_E =
u_1+u_2$ with
$$\gamma_R(u_1) + \gamma_C(u_2)=1.\leqno(6)$$
In particular, this implies that $\gamma_2(Id_E)=1$,
(where $\gamma_2(.)$ denotes the norm of factorization through Hilbert space,
see e.g. [P1, chapter 2] for more background) whence that $E$ is
isometric to $\ell^n_2$ $(n=\dim E)$. Moreover,
for any $e$ in the unit
sphere of $E$ we have
$$1 = \|e\|\le \|u_1(e)\| + \|u_2(e)\| \le \|u_1\| + \|u_2\| \le
\gamma_R(u_1) + \gamma_C(u_2)\le 1.$$
Therefore we must have
$$\|u_1(e)\| = \|u_1\| = \gamma_R(u_1) \hbox{ and } \|u_2(e)\| =
\|u_2\|=\gamma_C(u_2).\leqno{(7)}$$
Let $\alpha_i = \|u_i\|$ so that (by (6))
$\alpha_1+\alpha_2=1$. Assume that both
$\alpha_1>0$ and
$\alpha_2>0$. We will show that this is
impossible if $n>1$. Indeed, then
$U_i = (\alpha_i)^{-1}u_i$
$(i=1,2)$ is an isometry on
$\ell^n_2$, such that, for any $e$ in the unit
sphere of $E$, we have
$ e  =  \alpha_1U_1(e) + \alpha_2U_2(e) .$
By the strict convexity of $\ell^n_2$, this implies that $U_1(e) = U_2(e)
= e$ for all $e$. Moreover, by (7) we have
$\gamma_R(U_1) =1$ and
$\gamma_C(U_2)=1$. This implies that $E=R_n$ and $E=C_n$ completely
isometrically, which is absurd when $n>1$. Hence, if $n>1$, we conclude
that either $\alpha_1=0$ or $\alpha_2=0$, which
implies either
$\gamma_C(Id_E)=1$ or $\gamma_R(Id_E)=1$, equivalently either $E=C_n$ or
$E=R_n$ completely isometrically.
The remaining case $n=1$ is trivial.\qed

 \n{\bf Remark 11.} Alternate proof: if $\mu(E)=1$, then $E$ is
an injective operator space as well as its dual, but in [Ru2], Ruan gives
the complete list of the finite dimensional injective operator spaces.
Running down the list, and using an unpublished result of R. Smith
saying that a finite dimensional injective operator space is completely
contractively complemented in a finite dimensional $C^*$-algebra (see [B2]),
we find that $R_n$ and $C_n$ are the only
possibilities.

\n{\bf Remark 12.} We suspected that there   did not exist an operator space $X$
such that (with the notation of Corollary 9) we had for any $E_1,E_2$ and any $u\in E^*_1\otimes
E_2$
$$\gamma_{X}(\tilde u) =
\|u\|_\mu,\leqno(8)$$ and  indeed C. Le Merdy has kindly provided us with an argument, as follows.
Let $X$ be such a space. Let $E$ be an arbitrary finite dimensional subspace of
$X$ and let $v_E\in E^*\otimes X$   denote the tensor representing the inclusion map
$\tilde v_E \colon \ E \to X$.  Then, by (3) and (8),
$\gamma_2(\tilde  v_E)\le  \|v_E\|_\mu=\gamma_{X}(\tilde v_E)=1$. By a well known ultraproduct
argument (\cf e.g. [P1, p. 22]), this implies that $X$ is isometric to a Hilbert space. But
then, a variant of the  proof of Corollary 10 shows that we must have either $X=R$ or
$X=C$ completely isometrically, and this is absurd.

However, (8)  is true up to equivalence if we take for $X$ the direct
sum of $R$ and $C$, in any reasonable way.
For instance, it is easy to check that for any
$u\in E^*_1\otimes E_2$ ($\tilde u\colon \ E_1\to E_2$ being the associated
finite rank operator) we have
$${1\over 2 }\|u\|_\mu \le \gamma_{R\oplus_1 C}(\tilde u) \le \|u\|_\mu \leqno
 {  (9)}$$

\n {\bf Remark 13.}
It follows from Theorem~1 and the projectivity of Haagerup tensor
product that $\otimes_\mu$ is also projective, i.e, if $q_i : F_i \raw E_i$
($i=1,2$) are quotient maps, so is $q_1 \otimes q_2 : F_1 \otimes_\mu F_2
\raw E_1 \otimes_\mu E_2$. On the other hand, $\otimes_\mu$ is not
injective. To show this, consider the identity operator $i_n : R_n \cap C_n
\raw R_n \cap C_n$ and the natural (completely isometric) embedding $j_n : R_n \cap C_n
\hookrightarrow R_n \oplus_\infty C_n$. Remark 12 implies that
$$
|| j_n i_n ||_{(R_n + C_n) \otimes_\mu (R_n \oplus_\infty C_n)} \leq 2 .
$$
However, 
by [HP, p. 912] we have 
$\gamma_{\cal K}(i_n)\ge (1+\sqrt{n})/2,$
hence by Corollary 9, we have 
$$\| i_n \|_\mu \ge (1+\sqrt{n})/2.$$
This proves that the tensor product $\otimes_\mu$ is not injective.

\n{\bf Remark 14.} The examples in [P1, chapter 10] imply that there are
(infinite dimensional)
operator spaces $E$ such that $E^*\otimes_{\min} E=E^*\otimes_{\mu} E$
with equivalent norms, but $E$ is not completely  isomorphic to $R$ or $C$,
 and actually
 (as a Banach space) $E$ is not isomorphic to any Hilbert space.
Thus (in the isomorphic case) the second part of
 Corollary 10 does
not seem to extend to the infinite dimensional setting without assuming
some kind of approximation property.

We recall that any Hilbert space $H$ (resp. $K$)
can be equipped with a column (resp. row)
operator space, by identifying $H$  (resp. $K$)
with $H_c=B(\comp, H)$ (resp. with $K_r=B(K^*, \comp)$).
Any operator space of this form will be called a ``column space" (resp. a ``row space").

We will use the following result from [O]:

\proclaim Theorem 15. ([O]) Let $E$ be an operator
space such that $Id_E$ can be factorized
completely boundedly through the direct sum
$X=H_c \oplus_1 K_r$ of
 a column space and a row space, (i.e.\ there are
$c.b.$ maps
$u\colon \ E\to X$ and $v\colon \ X\to E$ such
that
$Id_E = vu$), then there are subspaces $E_1\subset
H_c$ and $E_2\subset K_r$ such that $E$ is
completely isomorphic to $E_1\oplus_1 E_2$. More
precisely, if we have
$\|u\|_{cb}\|v\|_{cb}\le c$ for some number $c$, then we can find a
complete isomorphism $T\colon \ E\to E_1\oplus_1 E_2$
 such that $\|T\|_{cb}\|T^{-1}\|_{cb}\le f(c)$ where $f\colon \ {\bf R}_+
\to {\bf R}_+$ is a certain function $f$.

\proclaim Theorem 16. The following properties of an
  an operator space $E$ are equivalent:
\item{\rm (i)} For any operator space $F$, we have
$F\otimes_{\min} E=F\otimes_{\mu} E$ isomorphically.
\item{\rm (ii)} $E$ is completely isomorphic to the direct sum
 of a row space and a column
space.

\pf The implication (i)~$\Rightarrow$~(ii) is easy and left to the reader.
 Conversely, assume
(i). Then, a routine argument shows that there is a constant $K$
such that for all $F$ and all $u$ in $F\otimes E$ we have $\|u\|_\mu \le
K\|u\|_{\rm min}$. Let $S\subset E$ be an arbitrary finite dimensional
subspace let $j_S\colon \ S\to E$ be the inclusion map, and let $\hat j_S
\in S^*\otimes E$ be the associated tensor. Then we have by (9)
$\sup\limits_S \gamma_{R\oplus_1 C} (j_S) = \sup\limits_S \|\hat j_S\|_\mu
\le K$. By a routine ultraproduct argument, this implies that the
identity of $E$ can be written as in
Theorem 15 with $c=K$, thus we conclude
that $E\simeq E_1\oplus_1 E_2$ where $E_1$ is a row space and $E_2$ a
column space. Note that we obtain an isomorphism $T\colon \ E\to E_1\oplus
_1 E_2$ such that $\|T\|_{cb} \|T^{-1}\|_{cb} \le f(K)$.
In particular, if $E$ is finite dimensional, we find $T$ such that
$$\|T\|_{cb} \|T^{-1}\|_{cb} \le f(\|i_E\|_\mu).$$\qed

We now turn to a result at the root of the present investigation. Let $E$
be an operator space and let $A$ be a unital operator algebra. For any $x$
in $E\otimes A$, we define
$$\delta(x) = \sup\|\sigma\cdot\pi(x)\|$$
where the supremum runs over all pairs $(\sigma,\pi)$ where $\sigma\colon
\ E\to B(H)$ is a complete contraction, $\pi\colon \ A\to B(H)$ a
morphism and moreover $\sigma$ and $\pi$ have commuting ranges.

Let $E\otimes_\delta A$ be the completion of $E\otimes A$ for this norm.
We may clearly also view $E \otimes_\delta A$ as an operator space using
the embedding $x\to \bigoplus\limits_{(\sigma,\pi)} (\sigma\cdot \pi)(x)$
where the direct sum runs over all pairs as above.

Then we may state.

\proclaim Theorem 17. Consider the linear mapping
$q\colon \ A\otimes E\otimes A\to E\otimes A$
defined by
$$q(a\otimes e\otimes b) = e\otimes (ab).$$
This mapping $q$ defines a complete metric surjection from $A \otimes_h E
\otimes_h
A$ onto $E\otimes_\delta A$.
More precisely, for any $n$ and any $x$ in $M_n(E\otimes A)$ with
$\|x\|_{M_n(E\otimes_\delta A)} <1$, there is $\tilde x$ in $M_n(A\otimes
E \otimes A)$ with $\|\tilde x\|_{M_n(A\otimes_h E \otimes_h A)} < 1$ such
that $I_{M_n} \otimes q(\tilde x) = x$.

\n {\bf Remark.} This statement is due to the second author [P4] (who is
indebted to C.~Le~Merdy for observing this useful reformulation).
A proof (somewhat different from the original one in [P4])  can be given following the lines
of  the
second  proof of Theorem~1, so we skip it.
 This result yields simpler proofs and extensions of several
 statements concerning nuclear $C^*$- algebras. See the final version of [P4] for more
details on this
topic.  
 
\n {\bf Acknowledgement.} We are very grateful to C. Le Merdy for letting us include
several useful remarks  and his  example showing the non-associativity of the $\mu$-tensor
product.

\bigskip\bigskip

\centerline {\bf References}
 \bigskip

 \item{[B1]} D. Blecher. Tensor products of
 operator spaces II.
  Canadian J. Math. 44 (1992) 75-90.

\item{[B2]} D. Blecher. The standard
 dual of an operator space.
  Pacific J. Math. 153 (1992) 15-30.

\item{[B3]} D. Blecher. A completely bounded
characterization of operator algebras. Math. Ann. 303 (1995)  227-239.

\item{[BRS]} D. Blecher, Z. J. Ruan  and  A. Sinclair.
A characterization of operator algebras. J. Funct. Anal.
89 (1990) 188-201.

\item{[BP1]} D. Blecher and V. Paulsen.  Tensor products of
operator spaces   J. Funct. Anal.  99 (1991) 262-292.

\item{[BS]} D. Blecher  and R. Smith. The dual of the Haagerup tensor
 product. Journal London Math. Soc. 45 (1992) 126-144.

\item{[CES]}  E. Christensen,  E. Effros  and  A. Sinclair.
Completely bounded multilinear maps and
$C^*$-algebraic cohomology. Invent. Math. 90 (1987)
279-296.

\item{[CS1]} E. Christensen and A. Sinclair.
Representations of completely bounded multilinear operators.
J. Funct. Anal. 72 (1987) 151-181.

\item{[CS2]}  E. Christensen  and A. Sinclair.  A survey of
completely bounded operators.  Bull. London Math. Soc.
21 (1989) 417-448.

\item{[Di]} P. Dixon.  Varieties of
 Banach algebras.
 Quarterly J. Math. Oxford 27 (1976) 481-487.

 \item{[ER1]} E. Effros  and Z.J. Ruan.  A new approach to
operators spaces.
 Canadian Math. Bull.
34 (1991) 329-337.

 \item{[ER2]} E. Effros and Z.J. Ruan.On the abstract
characterization of
 operator spaces.  Proc. Amer. Math. Soc. 119 (1993) 579-584.

\item{[ER3]} E. Effros and Z.J. Ruan. On approximation
properties for operator spaces, International J. Math. 1
(1990) 163-187.

 \item{[ER4]} E. Effros and Z.J. Ruan. Self duality for the Haagerup
 tensor product and Hilbert space factorization.  J. Funct. Anal. 100
(1991) 257-284.

 \item{[ER5]} E. Effros and Z.J. Ruan. Mapping spaces
and liftings for operator spaces.  Proc. London
Math. Soc.  69 (1994) 171-197.

\item{[ER6]}  E. Effros and Z.J. Ruan. The
Grothendieck-Pietsch and Dvoretzky-Rogers Theorems for
operator spaces. J. Funct. Anal. 122 (1994) 428-450.

\item{[HP]}  U. Haagerup and G. Pisier.  Bounded linear operators
between
$C^*-$algebras.
Duke Math. J. 71 (1993)  889-925.

 \item{[J]} M. Junge. Factorization theory
 for spaces of operators.
Habilitationsschrift, Kiel 1996.

 \item{[JP]} M. Junge and G. Pisier. Bilinear
 forms on exact operator
spaces and $B(H)\otimes B(H)$. Geometric and
Functional Analysis (GAFA Journal)
5 (1995) 329-363.

\item{[KaR]} R. Kadison  and J. Ringrose.   Fundamentals of
the theory
of operator algebras, vol. II   Academic Press.
New-York, 1986.

 \item{[Ki]} E. Kirchberg. On non-semisplit
extensions, tensor products and exactness of group
$C^*$-algebras. Invent. Math. 112 (1993) 449-489.

 \item{[O]} T. Oikhberg. Article in preparation.

\item{[Pa1]} V. Paulsen.   Completely bounded maps and
dilations.  Pitman Research Notes in
Math. 146, Longman, Wiley, New York, 1986.

\item{[Pa2]} V. Paulsen. Representation of
Function algebras, Abstract
 operator spaces and Banach space Geometry. J.
Funct. Anal. 109 (1992) 113-129.

\item{[Pa3]}  V. Paulsen.   The maximal operator space of a normed space.
  Proc. Edinburgh Math. Soc.  39 (1996) 309-323.

 \item{[PaP]}  V. Paulsen and S. Power.
Tensor products of non-self adjoint operator algebras.
 Rocky Mountain J.  Math. 20 (1990) 331-350.

 \item{[PaS]} V. Paulsen and R. Smith. Multilinear maps and
tensor norms on operator systems. J. Funct. Anal. 73
(1987) 258-276.

 \item{[Pes]} V. Pestov. Operator spaces and residually
finite-dimensional $C^*$-algebras.
J. Funct. Anal. 123 (1994) 308-317.

 \item{[P1]} G. Pisier.   Factorization of linear operators
and the Geometry of Banach spaces. CBMS (Regional
conferences of the A.M.S.)   60, (1986), Reprinted
with corrections 1987.

 \item{[P2]} G. Pisier. The operator Hilbert space $OH$,
complex interpolation and tensor norms.
Memoirs Amer.
Math. Soc.  vol. 122 , 585 (1996) 1-103.

\item{[P3]} G. Pisier.  A simple proof
of a theorem of Kirchberg
and related results on
$C^*$-norms.
J. Op. Theory.  35 (1996) 317-335.

\item{[P4]} G. Pisier. An introduction to the theory of
operator spaces.
Preprint. To appear.

\item{[Ru1]} Z.J. Ruan. Subspaces of $C^*$-algebras.
J. Funct. Anal. 76
 (1988) 217-230.

\item{[Ru2]} Z.J. Ruan. Injectivity and operator
spaces. Trans. Amer. Math. Soc. 315 (1989) 89-104.

 \item{[T]} M. Takesaki. Theory of Operator Algebras I.
Springer-Verlag New-York 1979.

\end
 values in non-commutative

\bye